\documentclass[a4paper,twoside]{article}
\usepackage{amsmath,amssymb,amsfonts}
\usepackage{graphics}
\newtheorem{proposition}{Proposition}[section]

\def\R{{\mathbb R}}

\date{\small\em \today}

\begin{document}
\title{Information geometric neighbourhoods of randomness and
geometry of the McKay bivariate gamma 3-manifold}
\author{ Khadiga Arwini and C.T.J. Dodson \\
 {\small\tt K.Alrawini@postgrad.umist.ac.uk \ \ \ \ dodson@umist.ac.uk}\\
{\small Department of Mathematics}\\
{\small University of Manchester Institute of Science and Technology} \\
{\small Manchester M60 1QD, UK}   }

\maketitle

\begin{abstract}
We show that gamma distributions provide models for departures
from randomness since every neighbourhood of an exponential
distribution contains a neighbourhood of gamma distributions,
using an information theoretic metric topology. We derive also the
information geometry of the 3-manifold of McKay bivariate gamma
distributions, which can provide a metrization of departures from
randomness and departures from
independence for bivariate processes. The curvature
objects are derived, including those on three submanifolds. As in
the case of bivariate normal manifolds, we have negative scalar
curvature but here it is not constant and we show how it depends
on correlation. These results have applications, for
example, in the characterization of stochastic materials.
\end{abstract}

\section{Gamma distributions and randomness}
The family of gamma probability density functions is given by
\begin{equation}
\{p(x;\beta,\alpha)= \left(\frac{\alpha}{\beta}\right)^\alpha \,
\frac{x^{\alpha-1}}{\Gamma(\alpha)}\, e^{-\frac{\alpha}{\beta}
\,x}| \alpha,\beta\in \R^+\}, \ \ x\in\R^+ \label{gamma}
\end{equation}
so the space of parameters is topologically $\R^+\times\R^+.$ It
is an exponential family and it includes as a special case
($\alpha=1$) the exponential distribution itself, which
complements the Poisson process on a line. It is pertinent to our
interests that the property of having sample standard deviation
independent of the mean actually characterizes gamma
distributions, as shown recently by Hwang and Hu~\cite{hwang}.
They proved, for $n\geq 3$ independent positive random variables
$x_1,x_2,\ldots , x_n$ with a common continuous probability
density function $f,$ that having independence of the sample mean
$\bar{x}$ and sample coefficient of variation $cv=S/\bar{x}$ is
equivalent to $f$ being a gamma distribution. Of course, the
exponential distribution has unit coefficient of variation.

\begin{figure}
\begin{picture}(300,220)(0,0)
\put(-20,0){\includegraphics{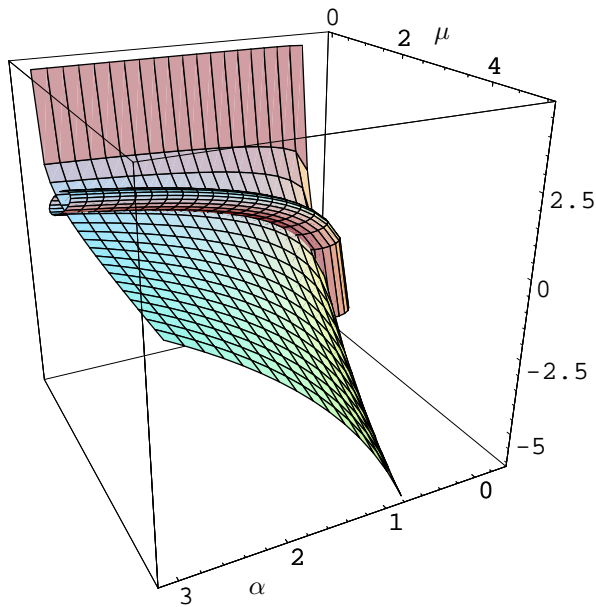}}
\put(250,220){$\mu$}
\put(180,10){$\alpha$}
\end{picture}
\caption{{\em Affine immersion in natural coordinates $\mu=\alpha/\beta,\alpha$
as a surface in $\R^3$ for the gamma manifold ${\cal G};$
the tubular neighbourhood surrounds all exponential distributions---these lie on
the curve $\alpha=1$ in the surface.}}
\label{expnhdnat}
\end{figure}

\begin{figure}
\begin{picture}(300,220)(0,0)
\put(-20,0){\includegraphics{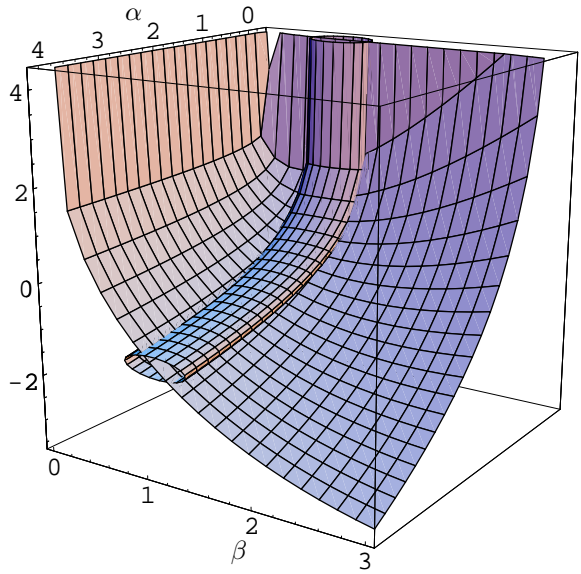}}
\put(140,214){$\alpha$}
\put(180,10){$\beta$}
\end{picture}
\caption{{\em Continuous image of the affine immersion
in Figure~\ref{expnhdnat} as a surface
in $\R^3$ using standard coordinates for the gamma manifold ${\cal G};$
the tubular neighbourhood surrounds all exponential distributions---these lie on
the curve $\alpha=1$ in the surface.}}
\label{expnhd}
\end{figure}

The univariate gamma distribution is widely used to model processes involving
a continuous positive random variable. Its information geometry is known and
has been applied recently to represent and metrize
departures from randomness of, for example,
the processes that allocate gaps between occurrences of each amino acid
along a protein chain within the {\em Saccharomyces cerevisiae} genome,
see Cai et al~\cite{amino}, clustering of galaxies and communications,
Dodson~\cite{ijtp,gsis,igcc}. In fact, we can make rather precise the
statement that around every random process on the real line there
is a neighbourhood of processes governed by the gamma distribution,
so gamma distributions can approximate any small enough
departure from randomness.
\begin{proposition}
Every neighbourhood of an exponential distribution contains a neighbourhood
of gamma distributions, using the subspace topology of $\R^3$ and information
theoretic immersions.
\end{proposition}
{\bf Proof:}
Dodson and Matsuzoe~\cite{affimm} have provided an affine immersion
in Euclidean $\R^3$ for ${\cal G},$ the manifold of gamma distributions with
Fisher information metric.
The coordinates $(\theta^1,\theta^2) = (\mu = \alpha/\beta, \alpha)$ form a
natural coordinate system (cf Amari and Nagaoka~\cite{AmNag}) for the gamma manifold ${\cal G}.$
Then ${\cal G}$ can be realized in Euclidean ${\R}^3$ by the graph of
the affine immersion $\{h,\xi\}$ where $\xi$ is the
 transversal vector field
along $h$~\cite{AmNag,affimm}:
\[
  h: {\cal G} \rightarrow \R^3 :
  \left( \! \!
    \begin{array}{c}
       \mu \\ \alpha
    \end{array} \! \! \right)
    \mapsto
    \left( \! \!
    \begin{array}{c}
       \mu \\ \alpha \\ \log\Gamma(\alpha) - \alpha\log\mu
    \end{array} \! \! \right), \quad
  \xi = \left( \! \! \begin{array}{c}
       0 \\ 0 \\ 1
    \end{array} \! \! \right).
\]
The submanifold of exponential distributions is represented by the curve
$$(0,\infty)\rightarrow \R^3: \mu \mapsto \{\mu,1,\log\frac{1}{\mu}\}$$
and a tubular neighbourhood of this curve will contain all immersions for small
enough perturbations of exponential distributions.
In Figure~\ref{expnhdnat} this is depicted in natural coordinates $\mu,\alpha$ and in
Figure~\ref{expnhd} the corresponding surface and tubular neighbourhood
(not here an affine immersion, just a continuous image) is shown in the usual
$(\alpha,\beta)$ coordinates of
the gamma family (\ref{gamma}).
The tubular neighbourhood in Figure~\ref{expnhd}
intersects with the gamma manifold immersion
to yield the required neighbourhood of gamma distributions, which completes our proof.
\hfill $\Box$

A simple transformation of random variable $x\in\R^+$ in (\ref{gamma}) to
$N=e^{-x}\in[0,1]$ converts a gamma distribution
into a log-gamma distribution, which turns out to have the same geometry.
\begin{proposition}[Dodson~\cite{gsis}]
The family of log-gamma probability density
functions
\begin{equation}
\{g(N,\alpha,\beta) = {\frac{{{{\frac{1}{N}}}^
       {1 - {\frac{\alpha }{\beta }}}}\,
     ({{{\frac{\alpha }{\beta }}})^{\alpha }}\,
     {{(\log {\frac{1}{N}})}^
       { \alpha -1}}}{
      \Gamma(\alpha )}}| \alpha,\beta\in \R^+\}, \ \ N\in[0,1]
  \label{loggamma}    \end{equation}
determines a Riemannian manifold $\cal{L}$ with information-theoretic metric
 having the properties:\\
$\bullet$ $\cal{L}$ contains the uniform distribution as the limit: \
$\lim_{\beta\rightarrow 1}g(N,\beta,1) = g(N,1,1) = 1$\\
$\bullet$ $\cal{L}$ contains approximations to truncated Normal
distributions for $\beta >>1$\\
$\bullet$ $\cal{L}$ is isometrically equivalent to the gamma manifold $\cal{G}$.
      \hfill $\Box$
\end{proposition}
Through this isometry and the result
of Dodson and Matsuzoe~\cite{affimm} we have an immersion in $\R^3$ that
represents also the log-gamma manifold and, since the isometry sends
the exponential distribution to the uniform distribution on $[0,1],$
we obtain another deduction:
\begin{proposition}
Every neighbourhood of the uniform distribution on $[0,1]$ contains
a neighbourhood of log-gamma distributions.            \hfill $\Box$
\end{proposition}
The value of such topological results lies in the fact that they
have qualitative consequences that are stable under small
perturbations of a process, something that would be important
in real applications. It gives confidence in the use of gamma
distributions to model near random processes. Moreover, the fact that
we have an information-theoretic metric topology for the neighbourhoods
lends significance to the result.

For comparison purposes, we recall that the
information geometry of the univariate and multivariate normal distributions
also are known, see for example Lauritzen~\cite{lauritzen} in Amari et al~\cite
{Am}
for the univariate case, Sato et al~\cite{sato} for the bivariate case
and Skovgaard~\cite{skov} for the multivariate case. We note in particular
that the univariate and bivariate normal distributions have constant negative
scalar curvature, so geometrically they constitute parts of pseudospheres.

First consider the 3-parameter family of
univariate gamma distributions with density function :
\begin{eqnarray}
  g(x;\beta,\alpha,\gamma) =\left(\frac{\alpha }{\beta } \right)^{\alpha }
\,\frac{\left( x - \gamma  \right) ^{ \alpha-1 }}
  {\Gamma(\alpha )}{e^{-\frac{\alpha \,\left( x - \gamma  \right) }{\beta }}},
\quad x>\gamma\geq0 ,\quad \beta,\alpha>0.
\label{threegamma}  \end{eqnarray} Evidently, the extra parameter
$\gamma \geq 0$ is a location shift and when $\gamma =0$ we
recover the univariate gamma distribution (\ref{gamma}), and when
$\alpha=1$ we obtain the exponential distributions with two
parameters.The mean $ \bar{x} $, standard deviation $\sigma_{x}$,
and coefficient of variation $c\nu_{x}$, for (\ref{threegamma})
are given by
$$ \bar{x}= \beta  + \gamma, \ \  \sigma_{x}^2 = \frac{{\beta }^2}{\alpha },
\ \ \ c\nu_{x}=\frac{\beta}{\sqrt{\alpha}(\beta + \gamma)}.
$$
The distribution (\ref{threegamma}) gives us a slight generalisation of the
gamma distribution (\ref{gamma}), which we use in the sequel.

One of the earliest forms of the bivariate gamma distribution is
due to Mckay~\cite{mckay}, defined by the density function
\begin{equation}
  f(x,y) =  \frac{c^{(\alpha_{1}+\alpha_{2})}x^{\alpha_{1}-1}(y-x)^{\alpha_{2}-
1}
e^{-c y}}{\Gamma(\alpha_{1})\Gamma(\alpha_{2})}\,\ {\rm defined\,
on} \,\ y>x>0 \ , \ \alpha_{1},c,\alpha_{2}>0
 \label{mckaybivariate}
 \end{equation}
The marginal distributions are gamma with shape parameters
$\alpha_{1}$ and $\alpha_{1}+\alpha_{2}$ , respectively.  The
covariance $Cov$ and correlation coefficient $\rho$ of $X$ and $Y$ are
given by :
\begin{eqnarray}
Cov(X,Y)&=&\frac{\alpha_{1}}{c^2} \\
\rho(X,Y)&=&\sqrt{\frac{\alpha_{1}}{\alpha_{1}+\alpha_{2}}}.
\end{eqnarray}

\section{Bivariate 3-parameter gamma 5-manifold}
\label{genMcKay} In this section we introduce a
bivariate gamma distribution which is a slight generalization of that
due to Mckay, by substituting $(x-\gamma_{1})$
for $x$ and $(y-\gamma_{2})$ for $y$ in equation
(\ref{mckaybivariate}). We call this a bivariate
3-parameter gamma distribution, because the marginal functions
are univariate 3-parameter gamma distributions.  Then we
consider the bivariate 3-parameter gamma models as a Riemannian
5-manifold.  The Christoffel symbols have been calculated but are
not listed because they are somewhat cumbersome.
\begin{proposition}
Let $X$ and $Y$ be continuous random variables, then
\begin{eqnarray}
  f(x,y) =  \frac{c^{(\alpha_{1}+\alpha_{2})}(x-\gamma_{1})^{\alpha_{1}-1}(y-
\gamma_{2}-x+\gamma_{1})^{\alpha_{2}-1}
e^{-c (y-\gamma_{2})}}{\Gamma(\alpha_{1})\Gamma(\alpha_{2})}\,\
 \label{arwinibivariate}  \end{eqnarray}
defined on $(y-\gamma_{2})>(x-\gamma_{1})>0 \ , \
\alpha_{1},c,\alpha_{2}>0,\,\gamma_{1},\gamma_{2}\geq0,$ is a
density function.  The covariance and
marginal density functions, of $X$ and $Y$ are given by:
\begin{eqnarray}
 \sigma_{12}&=&\frac{\alpha_{1}}{c^{2}}  \\
 f_{X}(x) &= & \frac{c^{\alpha_{1}} (x-\gamma_{1})^
{\alpha_{1}-1}e^{-c (x-\gamma_{1})}}{\Gamma(\alpha_{1})},\quad x>\gamma_{1}
\geq0 \\
 f_{Y}(y)&=&\frac{c^{(\alpha_{1}+\alpha_{2})}
(y-\gamma_{2})^{(\alpha_{1}+\alpha_{2})-1}e^
{-c (y-\gamma_{2})}}{\Gamma(\alpha_{1}+\alpha_{2})},
\quad y>\gamma_{2}\geq0
\end{eqnarray}
$\hfill \Box$
\end{proposition}
Note that the marginal functions $ f_{X}$ and $ f_{Y}$ are
univariate 3-parameter gamma distributions with parameters
$(c,\alpha_{1},\gamma_{1})$ and
$(c,\alpha_{1}+\alpha_{2},\gamma_{2})$, where $\gamma_{1} $ and
$\gamma_{2} $ are location parameters.
We shall refer to (\ref{arwinibivariate}) as giving the
bivariate 3-parameter gamma distributions.

It is easily shown that $ c
=\sqrt{\frac{\alpha_{1}}{\sigma_{12}}} $ , so the bivariate
3-parameter gamma distribution (\ref{arwinibivariate}) can be
presented in the form:
\begin{eqnarray}
f(x,y)={\left( \frac{{{\alpha }_1}}{{{\sigma }_{12}}} \right)
}^{\frac{{{\alpha }_1} + {{\alpha }_2}}{2}}\, \frac{{\left( x -
{{\gamma }_1} \right) }^{ {{\alpha }_1}-1}\, {\left(  y  -
{{\gamma }_2}+x- {{\gamma }_1} \right) }^{ {{\alpha }_2}-1}\,}{
    \Gamma({{\alpha }_1})\,\Gamma({{\alpha
    }_2})}\,e^{-{\sqrt{\frac{{{\alpha }_1}}{{{\sigma }_{12}}}}} \,\left( y -
{{\gamma }_2} \right)},
    \label{arwinidistribution}
\end{eqnarray}
defined on $(y-\gamma_{2})>(x-\gamma_{1})>0,$\, with parameters
$\alpha_{1},\alpha_{2},\sigma_{12}>0,$\, and
$\gamma_{1},\gamma_{2} \geq0.$

\begin{proposition}
Let $M^{*} $ be the set of bivariate 3-parameter gamma
distributions, that is
\begin{eqnarray}
M^{*}& =&  \{f |f(x,y) ={\left( \frac{{{\alpha }_1}}{{{\sigma
}_{12}}} \right) }^{\frac{{{\alpha }_1} + {{\alpha }_2}}{2}}
\frac{{\left( x - {{\gamma }_1} \right) }^{ {{\alpha }_1}-1}\,
{\left(  y  - {{\gamma }_2}+x- {{\gamma }_1} \right) }^{ {{\alpha
}_2}-1}\,}{
    \Gamma({{\alpha }_1})\,\Gamma({{\alpha
    }_2})}\,e^{-{\sqrt{\frac{{{\alpha }_1}}{{{\sigma }_{12}}}}} \,\left( y -
{{\gamma }_2}
    \right)},\nonumber  \\
   && \quad (y-\gamma_{2})>(x-\gamma_{1})>0,\,
\alpha_{1},\alpha_{2}>2,\,\sigma_{12}>0,\,\gamma_{1},\gamma_{2}
\geq0.\} \label{arwinimodel} \end{eqnarray}
 Then we have :
\begin{enumerate}
\item Identifying $(\alpha_{1},\alpha_{2},\sigma_{12},\gamma_{1},\gamma_{2})$
as a local coordinate
system, $M^{*}$ can be regarded as a 5-manifold.
\item $M^{*}$ is a Riemannian manifold with Fisher information matrix $
G=[g_{ij}]$ where
\begin{displaymath}
g_{ij}=
\int_{\gamma_{1}}^{\infty}\int_{x-\gamma_{1}+\gamma_{2}}^{\infty}\frac{\partial^
{2}\log
f(x,y;x^{1},x^{2},x^{3},x^{4},x^{5})}{\partial x^{i}\partial
x^{j}}\ f(x,y;x^{1},x^{2},x^{3},x^{4},x^{5}) \ dy\  dx
\end{displaymath}
and  $ x^{1}=\alpha_{1} ,\,x^{2}=\alpha_{2},\,
x^{3}=\sigma_{12},\,
x^{4}=\gamma_{1} ,\,x^{5}=\gamma_{2} .$\\
is given by :
\begin{eqnarray}
 \left[\begin{array}{ccccc} \psi'({\alpha}_1) +
\frac{-3\,{{{\alpha }}_1} +
     {{{\alpha }}_2}}{4\,{{{\alpha }_1}}^2}  &
\frac{-1}{2\,{{\alpha }_1}}  & \frac{{{\alpha }_1} - {{\alpha
}_2}}{4\,{{\alpha }_1}\,{{\sigma }_{12}}} &
 \frac{{\sqrt{{{\alpha }_1}}}}{\left( -1 + {{\alpha }_1} \right) \,{\sqrt
{{{\sigma }_{12}}}}} &
 \frac{-1}{2\,{\sqrt{{{\alpha }_1}}}\,{\sqrt{{{\sigma }_{12}}}}} \\
\frac{-1}{2\,{{\alpha }_1}} &
 \psi'({\alpha}_2) &
 \frac{1}{2\,{{\sigma }_{12}}}  &
  \frac{{\sqrt{{{\alpha }_1}}}}{\left( 1 - {{\alpha }_2} \right) \,{\sqrt
{{{\sigma }_{12}}}}}  &
 \frac{{\sqrt{{{\alpha }_1}}}}{\left( -1 + {{\alpha }_2} \right) \,{\sqrt
{{{\sigma }_{12}}}}}  \\
\frac{{{\alpha}_1}-{{\alpha}_2}}{4\,{{\alpha}_1}\,{{\sigma}_{12}}}
 &
 \frac{1}{2\,{{\sigma }_{12}}}  &
 \frac{{{\alpha }_1} + {{\alpha }_2}}{4\,{{{\sigma }_{12}}}^2} & 0  &
 \frac{{\sqrt{{{\alpha }_1}}}}{2\,{{{\sigma }_{12}}}^{\frac{3}{2}}} \\
 \frac{{\sqrt{{{\alpha }_1}}}}{\left( -1 + {{\alpha }_1} \right) \,{\sqrt
{{{\sigma }_{12}}}}} &
\frac{{\sqrt{{{\alpha }_1}}}}{\left( 1 - {{\alpha }_2} \right)
\,{\sqrt{{{\sigma }_{12}}}}}  & 0 &
 \frac{{{\alpha }_1}}{\left( -2 + {{\alpha }_1} \right)
\,{{\sigma }_{12}}} +
  \frac{{{\alpha }_1}}{\left( -2 + {{\alpha }_2} \right) \,{{\sigma }_{12}}} &
 \frac{-{{\alpha }_1}}{\left( -2 + {{\alpha }_2} \right) \,{{\sigma
}_{12}}}  \\
\frac{-1}{2\,{\sqrt{{{\alpha }_1}}}\,{\sqrt{{{\sigma }_{12}}}}} &
\frac{{\sqrt{{{\alpha }_1}}}}{\left( -1 +
{{\alpha}_2}\right)\,{\sqrt{{{\sigma }_{12}}}}}  &
 \frac{{\sqrt{{{\alpha }_1}}}}{2\,{{{\sigma
 }_{12}}}^{\frac{3}{2}}} &
\frac{-{{\alpha }_1}}{\left( -2 + {{\alpha }_2} \right) \,{{\sigma
}_{12}}}  &
  \frac{{{\alpha }_1}}{\left( -2 + {{\alpha }_2} \right) \,{{\sigma }_{12}}}
\end{array}\right]
\label{arwinimetric}
\end{eqnarray}
where $
\psi(\alpha_{i})=\frac{\Gamma'(\alpha_{i})}{\Gamma(\alpha_{i})}\quad(i=1,2).
\hfill \Box $
\end{enumerate}
\end{proposition}

\section{Mckay bivariate gamma 3-manifold}\label{McKay}
In this section we consider the Mckay bivariate gamma model as a
3-manifold, equipped with Fisher information as Riemannian metric,
and derive the induced geometry, i.e., the Ricci tensor, the
scalar curvatures etc; the Christoffel symbols were computed but
are omitted here.  In addition, we consider three submanifolds as
special cases, and discuss their geometrical structure.

\subsubsection{Fisher information metric}
The classical family of Mckay bivariate gamma distributions is given by:
\begin{eqnarray}
f(x,y;\alpha_{1},\sigma_{12},\alpha_{2}) =
\frac{(\frac{\alpha_{1}}{\sigma_{12}})^{\frac{(\alpha_{1}+\alpha_{2})}{2}}x^
{\alpha_{1}-1}(y-x)^{\alpha_{2}-1}
e^{-\sqrt{\frac{\alpha_{1}}{\sigma_{12}}}y}}{\Gamma(\alpha_{1})\Gamma(\alpha_
{2})}
\ , \label{mckaydistribution}
\end{eqnarray}
defined on $ 0<x<y<\infty $ with parameters $
\alpha_{1},\sigma_{12},\alpha_{2}>0$.  Where $\sigma_{12}$ is the
covariance of $X$ and $Y$.
The correlation coefficient and
marginal functions, of $X$ and $Y$ are given by :
\begin{eqnarray}
  \rho(X,Y)&=& \sqrt{\frac{\alpha_{1}}{\alpha_{1}+\alpha_{2}}}\\
 f_{X}(x) &= & \frac{(\frac{\alpha_{1}}{\sigma_{12}})^{\frac{\alpha_{1}}{2}} x^
{\alpha_{1}-1}e^{-\sqrt{\frac{\alpha_{1}}{\sigma_{12}}}
x}}{\Gamma(\alpha_{1})},\quad x
>0 \\
 f_{Y}(y)&=&\frac{(\frac{\alpha_{1}}{\sigma_{12}})^{\frac{(\alpha_{1}+\alpha_{2})}{2}}
y^{(\alpha_{1}+\alpha_{2})-1}e^
{-\sqrt{\frac{\alpha_{1}}{\sigma_{12}}}
y}}{\Gamma(\alpha_{1}+\alpha_{2})}, \quad y>0
\end{eqnarray}
Note that it is not possible to choose parameters such that both
marginal functions are exponential.

\begin{proposition}
Let $M$ be the set of Mckay bivariate gamma distributions, that is
\begin{eqnarray}
M&=&\{f |f(x,y;\alpha_{1},\sigma_{12},\alpha_{2}) =
\frac{(\frac{\alpha_{1}}{\sigma_{12}})^{\frac{(\alpha_{1}+\alpha_{2})}{2}}x^
{\alpha_{1}-1}(y-x)^{\alpha_{2}-1}
e^{-\sqrt{\frac{\alpha_{1}}{\sigma_{12}}}y}}{\Gamma(\alpha_{1})\Gamma(\alpha_
{2})}\,,
\nonumber \\ &&  y>x>0,\ \alpha_{1},\sigma_{12},\alpha_{2}>0\}
  \label{mckaymodel} \end{eqnarray}
 Then we have :
\begin{enumerate}
\item Identifying $(\alpha_{1},\sigma_{12},\alpha_{2})$ as a local coordinate
system, $M$ is a 3-manifold.
\item  $M$ is a Riemannian 3-manifold with
Fisher information metric $G=[g_{ij}]$ given by :
\begin{eqnarray}
   [g_{ij}]=\left[ \begin{array}{ccc}
\frac{-3\,{\alpha_{1}} + {\alpha_{2}}}{4\,{\alpha_{1}}^2} + \psi'(\alpha_{1}) &
\frac{{\alpha_{1}} - {\alpha_{2}}}{4\,{\alpha_{1}}\,{\sigma_{12}}} & -\frac{1}{2
\,{\alpha_{1}}} \\
\frac{{\alpha_{1}} - {\alpha_{2}}}{4\,{\alpha_{1}}\,{{\sigma_
{12}}}}              & \frac{{\alpha_{1}}+ {\alpha_{2}}}{4\,{{\sigma_{12}}}
^2}       & \frac{1}{2\,{\sigma_{12}}}  \\
                               -\frac{1}{2\,{\alpha_{1}}}   & \frac{1}{2\,
{\sigma_{12}}}                      & \psi'(\alpha_{2})
\end{array} \right] \label{mckaymetric}
\end{eqnarray}
\item  The inverse $[g^{ij}]$ of  $[g_{ij}]$ is given by:
\begin{eqnarray}
  g^{11}&=& -\left( \frac{-1 + \left( {\alpha_{1}} +
{\alpha_{2}} \right) \,
       \psi'(\alpha_{2})}{\psi'(\alpha_{2}) +
      \psi'(\alpha_{1})\,
       \left( 1 - \left( {\alpha_{1}} + {\alpha_{2}} \right) \,
          \psi'(\alpha_{2}) \right) } \right)\, ,\nonumber\\
  g^{12}&=& g^{21}= \frac{{\sigma_{12}}\,\left( 1 + \left(
{\alpha_{1}} - {\alpha_{2}}\right) \,
       \psi'(\alpha_{2}) \right) }{{\alpha_{1}}\,
    \left( \psi'(\alpha_{2}) +
      \psi'(\alpha_{1})\,
       \left( 1 - \left( {\alpha_{1}} + {\alpha_{2}} \right) \,
          \psi'(\alpha_{2}) \right)  \right) }\, ,\nonumber \\
 g^{13}&=& g^{31}=\frac{1}{-\psi'(\alpha_{2}) +
    \psi'(\alpha_{1})\,
     \left( -1 + \left( {\alpha_{1}} + {\alpha_{2}} \right) \,
        \psi'(\alpha_{2}) \right) } \,,\nonumber\\
 g^{22}&=& \frac{{{\sigma_{12}}}^2\,\left( -1 + \left( -3\,{\alpha_{1}} +
{\alpha_{2}} +
         4\,{{\alpha_{1}}}^2\,\psi'(\alpha_{1}) \right) \,
       \psi'(\alpha_{2}) \right) }{{{\alpha_{1}}}^2\,
    \left( -\psi'(\alpha_{2}) +
      \psi'(\alpha_{1})\,
       \left( -1 + \left( {\alpha_{1}} + {\alpha_{2}} \right) \,
         \psi'(\alpha_{2}) \right)  \right) } \,,\nonumber\\
  g^{23}&=& g^{32}=\frac{{\sigma_{12}}\,\left( -1 +
2\,{\alpha_{1}}\,\psi'(\alpha_{1}) \right) }
  {{\alpha_{1}}\,\left( \psi'(\alpha_{2}) +
      \psi'(\alpha_{1})\,
       \left( 1 - \left( {\alpha_{1}} + {\alpha_{2}} \right) \,
          \psi'(\alpha_{2}) \right)  \right) }\, ,\nonumber\\
 g^{33}&=&-\left( \frac{-1 + \left( {\alpha_{1}} + {\alpha_{2}}
\right) \,
       \psi'(\alpha_{1})}{\psi'(\alpha_{2}) +
      \psi'(\alpha_{1})\,
       \left( 1 - \left( {\alpha_{1}} + {\alpha_{2}} \right) \,
          \psi'(\alpha_{2}) \right) } \right)\,.
          \end{eqnarray}
\end{enumerate}
$\hfill \Box$
\end{proposition}

\subsubsection{Curvature properties}
We provide the various curvature objects of the McKay 3-manifold $M$;
the Christoffel symbols are known but they are omitted here.

\begin{proposition}The components of the curvature tensor $ R_{ijkl}$ are given
by:
\begin{eqnarray}
  R_{1212}&=&
\frac{\psi'({\alpha_{2}})\,\left( \psi'({\alpha_{1}}) +
      \left( {\alpha_{1}} + {\alpha_{2}} \right) \,\psi''({\alpha_{1}})
\right) }{16\,{{\sigma_{12}}}^2\,
    \left( \psi'({\alpha_{1}}) + \psi'({\alpha_{2}}) -
      \left( {\alpha_{1}} + {\alpha_{2}} \right) \,\psi'({\alpha_{1}})
\,\psi'({\alpha_{2}})
      \right) }  \, ,\nonumber\\
 R_{1213}&=& \frac{\psi'({\alpha_{2}})\,\left( \psi'({\alpha_{1}}) +
      2\,{\alpha_{1}}\,\psi''({\alpha_{1}}) \right) }{16\,{\alpha_{1}}\,{\sigma_
{12}}\,
    \left( \psi'({\alpha_{1}}) + \psi'({\alpha_{2}}) -
      \left( {\alpha_{1}} + {\alpha_{2}} \right) \,\psi'({\alpha_{1}})
\,\psi'({\alpha_{2}})
      \right) }  \, ,\nonumber\\
 R_{1223}&=& \frac{- \psi'({\alpha_{1}})\,\psi'({\alpha_{2}}) }
  {16\,{{\sigma_{12}}}^2\,\left( \psi'({\alpha_{1}}) + \psi'({\alpha_{2}}) -
      \left( {\alpha_{1}} + {\alpha_{2}} \right) \,\psi'({\alpha_{1}})
\,\psi'({\alpha_{2}})
      \right) } \,  ,\nonumber\\
 R_{1313}&=& \frac{-\left( -\left(
\psi'({\alpha_{1}})\,\psi'({\alpha_{2}}) \right)  +
      \left( \left( 3\,{\alpha_{1}} - {\alpha_{2}} \right) \,\psi'({\alpha_
{1}}) +
         4\,{{\alpha_{1}}}^2\,\psi''({\alpha_{1}}) \right) \,\psi''({\alpha_
{2}}) \right) }{16\,
    {{\alpha_{1}}}^2\,\left( \psi'({\alpha_{1}}) + \psi'({\alpha_{2}}) -
      \left( {\alpha_{1}} + {\alpha_{2}} \right) \,\psi'({\alpha_{1}})
\,\psi'({\alpha_{2}})
      \right) }  \, ,\nonumber\\
 R_{1323}&=& \frac{- \psi'({\alpha_{2}})\,
      \left( \psi'({\alpha_{2}}) +
        \left( -{\alpha_{1}} + {\alpha_{2}} \right) \,\psi''({\alpha_{2}})
\right)
      }{16\,{\alpha_{1}}\,{\sigma_{12}}\,\left( \psi'({\alpha_{1}}) +
\psi'({\alpha_{2}}) -
      \left( {\alpha_{1}} + {\alpha_{2}} \right) \,\psi'({\alpha_{1}})
\,\psi'({\alpha_{2}})
      \right) } \,  ,\nonumber\\
 R_{2323}&=& \frac{\psi'({\alpha_{1}})\,\left( \psi'({\alpha_{2}})+
      \left( {\alpha_{1}} + {\alpha_{2}} \right) \,\psi''({\alpha_{2}})
\right) }{16\,{{\sigma_{12}}}^2\,
    \left( \psi'({\alpha_{1}}) + \psi'({\alpha_{2}}) -
      \left( {\alpha_{1}} + {\alpha_{2}} \right) \,\psi'({\alpha_{1}})
\,\psi'({\alpha_{2}})
      \right) }   \,,
      \end{eqnarray}
      while the other independent components are zero.
$\hfill \Box$
\end{proposition}

\begin{proposition} The components of the Ricci tensor are given
by the symmetric matrix $ R=[R_{ij}]$:
\begin{eqnarray*}
R_{11} &=&\frac{-3\,{\psi'({\alpha_{1}})}^2\,\psi'({\alpha_{2}}) -
     3\,\psi'({\alpha_{1}})\,{\psi'({\alpha_{2}})}^2}{16\,{\alpha_{1}}\,
     {\left( \psi'({\alpha_{1}}) + \psi'({\alpha_{2}}) -
         \left( {\alpha_{1}} + {\alpha_{2}} \right) \,\psi'({\alpha_{1}})\,
          \psi'({\alpha_{2}}) \right) }^2} + \nonumber\\&&
  \frac{{\alpha_{2}}\,{\psi'({\alpha_{1}})}^2\,\psi'({\alpha_{2}}) +
     {\alpha_{2}}\,\psi'({\alpha_{1}})\,{\psi'({\alpha_{2}})}^2}{16\,{{\alpha_
{1}}}^2\,
     {\left( \psi'({\alpha_{1}}) + \psi'({\alpha_{2}}) -
         \left( {\alpha_{1}} + {\alpha_{2}} \right) \,\psi'({\alpha_{1}})\,
          \psi'({\alpha_{2}}) \right) }^2} + \nonumber\\&&
  \frac{{\psi'({\alpha_{1}})}^2\,{\psi'({\alpha_{2}})}^2 -
     2\,\psi'({\alpha_{1}})\,\psi'({\alpha_{2}})\,\psi''({\alpha_{1}})}{4\,
     {\left( \psi'({\alpha_{1}}) + \psi'({\alpha_{2}}) -
         \left( {\alpha_{1}} + {\alpha_{2}} \right) \,\psi'({\alpha_{1}})\,
          \psi'({\alpha_{2}}) \right) }^2} + \nonumber\\&&
  \frac{-\left( {\alpha_{2}}\,\psi'({\alpha_{2}})\,\psi''({\alpha_{1}})
\right)  +
     {{\alpha_{2}}}^2\,{\psi'({\alpha_{2}})}^2\,\psi''({\alpha_{1}})}{16\,
{{\alpha_{1}}}^2\,
     {\left( \psi'({\alpha_{1}}) + \psi'({\alpha_{2}}) -
         \left( {\alpha_{1}} + {\alpha_{2}} \right) \,\psi'({\alpha_{1}})\,
          \psi'({\alpha_{2}}) \right) }^2} + \nonumber\\&&
  \frac{{\alpha_{1}}\,\psi'({\alpha_{1}})\,{\psi'({\alpha_{2}})}^2\,
      \psi''({\alpha_{1}}) + {\alpha_{2}}\,\psi'({\alpha_{1}})\,{\psi'({\alpha_
{2}})}^2\,
      \psi''({\alpha_{1}})}{4\,{\left( \psi'({\alpha_{1}}) + \psi'({\alpha_
{2}}) -
         \left( {\alpha_{1}} + {\alpha_{2}} \right) \,\psi'({\alpha_{1}})\,
          \psi'({\alpha_{2}}) \right) }^2} + \nonumber\\&&
  \frac{3\,\psi'({\alpha_{2}})\,\psi''({\alpha_{1}}) +
     3\,\psi'({\alpha_{1}})\,\psi''({\alpha_{2}})}{16\,{\alpha_{1}}\,
     {\left( \psi'({\alpha_{1}}) + \psi'({\alpha_{2}}) -
         \left( {\alpha_{1}} + {\alpha_{2}} \right) \,\psi'({\alpha_{1}})\,
          \psi'({\alpha_{2}}) \right) }^2} + \nonumber\\&&
  \frac{-3\,{\psi'({\alpha_{2}})}^2\,\psi''({\alpha_{1}}) -
     3\,{\psi'({\alpha_{1}})}^2\,\psi''({\alpha_{2}})}{16\,
     {\left( \psi'({\alpha_{1}}) + \psi'({\alpha_{2}}) -
         \left( {\alpha_{1}} + {\alpha_{2}} \right) \,\psi'({\alpha_{1}})\,
          \psi'({\alpha_{2}}) \right) }^2} + \nonumber\\&&
  \frac{-\left( {\alpha_{2}}\,{\psi'({\alpha_{2}})}^2\,\psi''({\alpha_{1}})
\right)  -
     {\alpha_{2}}\,{\psi'({\alpha_{1}})}^2\,\psi''({\alpha_{2}})}{8\,{\alpha_
{1}}\,
     {\left( \psi'({\alpha_{1}}) + \psi'({\alpha_{2}}) -
         \left( {\alpha_{1}} + {\alpha_{2}} \right) \,\psi'({\alpha_{1}})\,
          \psi'({\alpha_{2}}) \right) }^2} + \nonumber\\&&
  \frac{-\left( {\alpha_{2}}\,\psi'({\alpha_{1}})\,\psi''({\alpha_{2}})
\right)  +
     {{\alpha_{2}}}^2\,{\psi'({\alpha_{1}})}^2\,\psi''({\alpha_{2}})}{16\,
{{\alpha_{1}}}^2\,
     {\left( \psi'({\alpha_{1}}) + \psi'({\alpha_{2}}) -
         \left( {\alpha_{1}} + {\alpha_{2}} \right) \,\psi'({\alpha_{1}})\,
          \psi'({\alpha_{2}}) \right) }^2} + \nonumber\\&&
  \frac{\psi''({\alpha_{1}})\,\psi''({\alpha_{2}}) -
     {\alpha_{1}}\,\psi'({\alpha_{1}})\,\psi''({\alpha_{1}})\,\psi''({\alpha_
{2}}) -
     {\alpha_{2}}\,\psi'({\alpha_{1}})\,\psi''({\alpha_{1}})\,\psi''({\alpha_
{2}})}{4\,
     {\left( \psi'({\alpha_{1}}) + \psi'({\alpha_{2}}) -
         \left( {\alpha_{1}} + {\alpha_{2}} \right) \,\psi'({\alpha_{1}})\,
          \psi'({\alpha_{2}}) \right) }^2} \,,\nonumber
 \end{eqnarray*}
\begin{eqnarray*}
R_{12}&=& \frac{{\psi'({\alpha_{1}})}^2\,\psi'({\alpha_{2}}) +
     \psi'({\alpha_{1}})\,{\psi'({\alpha_{2}})}^2 -
     \psi'({\alpha_{2}})\,\psi''({\alpha_{1}})}{16\,{\sigma_{12}}\,
     {\left( \psi'({\alpha_{1}}) + \psi'({\alpha_{2}}) -
         \left( {\alpha_{1}} + {\alpha_{2}} \right) \,\psi'({\alpha_{1}})\,
          \psi'({\alpha_{2}}) \right) }^2} + \nonumber\\&&
  \frac{-\left( {\alpha_{2}}\,{\psi'({\alpha_{1}})}^2\,\psi'({\alpha_{2}})
\right)  -
     {\alpha_{2}}\,\psi'({\alpha_{1}})\,{\psi'({\alpha_{2}})}^2 +
     {\alpha_{2}}\,\psi'({\alpha_{2}})\,\psi''({\alpha_{1}})}{16\,{\alpha_{1}}\,
{\sigma_{12}}\,
     {\left( \psi'({\alpha_{1}}) + \psi'({\alpha_{2}}) -
         \left( {\alpha_{1}} + {\alpha_{2}} \right) \,\psi'({\alpha_{1}})\,
          \psi'({\alpha_{2}}) \right) }^2} +\nonumber\\&&
  \frac{{\alpha_{1}}\,{\psi'({\alpha_{2}})}^2\,\psi''({\alpha_{1}}) -
     \psi'({\alpha_{1}})\,\psi''({\alpha_{2}}) +
     {\alpha_{1}}\,{\psi'({\alpha_{1}})}^2\,\psi''({\alpha_{2}})}{16\,{\sigma_
{12}}\,
     {\left( \psi'({\alpha_{1}}) + \psi'({\alpha_{2}}) -
         \left( {\alpha_{1}} + {\alpha_{2}} \right) \,\psi'({\alpha_{1}})\,
          \psi'({\alpha_{2}}) \right) }^2} + \nonumber\\&&
  \frac{-\left( {{\alpha_{2}}}^2\,{\psi'({\alpha_{2}})}^2\,\psi''({\alpha_{1}})
\right)  +
     {\alpha_{2}}\,\psi'({\alpha_{1}})\,\psi''({\alpha_{2}}) -
     {{\alpha_{2}}}^2\,{\psi'({\alpha_{1}})}^2\,\psi''({\alpha_{2}})}{16\,
{\alpha_{1}}\,{\sigma_{12}}\,
     {\left( \psi'({\alpha_{1}}) + \psi'({\alpha_{2}}) -
         \left( {\alpha_{1}} + {\alpha_{2}} \right) \,\psi'({\alpha_{1}})\,
          \psi'({\alpha_{2}}) \right) }^2}\, ,\nonumber
 \end{eqnarray*}
\begin{eqnarray*}
 R_{13}&=& \frac{-\left( {\psi'({\alpha_{1}})}^2\,\psi'({\alpha_{2}}) \right)  -
     \psi'({\alpha_{1}})\,{\psi'({\alpha_{2}})}^2 +
     \psi'({\alpha_{2}})\,\psi''({\alpha_{1}})}{8\,{\alpha_{1}}\,
     {\left( \psi'({\alpha_{1}}) + \psi'({\alpha_{2}}) -
         \left( {\alpha_{1}} + {\alpha_{2}} \right) \,\psi'({\alpha_{1}})\,
          \psi'({\alpha_{2}}) \right) }^2} + \nonumber\\&&
  \frac{{\psi'({\alpha_{2}})}^2\, \left( \psi''({\alpha_{1}}) + 2\,
{\psi'({\alpha_{1}})}^2 \right)+
    \psi''({\alpha_{2}}) \left( {\psi'({\alpha_{1}})}^2\,+ 2\,\psi''({\alpha_
{1}})
    \right)}{8\,
     {\left( \psi'({\alpha_{1}}) + \psi'({\alpha_{2}}) -
         \left( {\alpha_{1}} + {\alpha_{2}} \right) \,\psi'({\alpha_{1}})\,
          \psi'({\alpha_{2}}) \right) }^2} +\nonumber\\&&
  \frac{-\left( {\alpha_{2}}\,{\psi'({\alpha_{2}})}^2\,\psi''({\alpha_{1}})
\right)  +
     \psi'({\alpha_{1}})\,\psi''({\alpha_{2}}) -
     {\alpha_{2}}\,{\psi'({\alpha_{1}})}^2\,\psi''({\alpha_{2}})}{8\,{\alpha_
{1}}\,
     {\left( \psi'({\alpha_{1}}) + \psi'({\alpha_{2}}) -
         \left( {\alpha_{1}} + {\alpha_{2}} \right) \,\psi'({\alpha_{1}})\,
          \psi'({\alpha_{2}}) \right) }^2} \,,\nonumber
\end{eqnarray*}
\begin{eqnarray*}
         R_{22}&=& \frac{\left( {\alpha_{1}} + {\alpha_{2}} \right)
\,\left( \psi'({\alpha_{2}})\,
       \left( \psi'({\alpha_{1}})\,
          \left( \psi'({\alpha_{1}}) + \psi'({\alpha_{2}}) \right)  -
         \psi''({\alpha_{1}}) + \left( {\alpha_{1}} + {\alpha_{2}} \right) \,
          \psi'({\alpha_{2}})\,\psi''({\alpha_{1}}) \right) \right) }
       {16\,{{\sigma_{12}}}^2\,
    {\left( \psi'({\alpha_{1}}) + \psi'({\alpha_{2}}) -
        \left( {\alpha_{1}} + {\alpha_{2}} \right) \,\psi'({\alpha_{1}})
\,\psi'({\alpha_{2}})
        \right) }^2} + \nonumber\\&&
        \frac{\left( {\alpha_{1}} + {\alpha_{2}} \right)
\,\left( \psi'({\alpha_{1}})\,\left( -1 +
         \left( {\alpha_{1}} + {\alpha_{2}} \right) \,\psi'({\alpha_{1}})
\right) \,
       \psi''({\alpha_{2}}) \right) } {16\,{{\sigma_{12}}}^2\,
    {\left( \psi'({\alpha_{1}}) + \psi'({\alpha_{2}}) -
        \left( {\alpha_{1}} + {\alpha_{2}} \right) \,\psi'({\alpha_{1}})
\,\psi'({\alpha_{2}})
        \right) }^2} \, , \nonumber
\end{eqnarray*}
\begin{eqnarray*}
  R_{23}&=& \frac{\psi'({\alpha_{2}})\,\left( \psi'({\alpha_{1}})\,
        \left( \psi'({\alpha_{1}}) + \psi'({\alpha_{2}}) \right)  -
       \psi''({\alpha_{1}}) + \left( {\alpha_{1}} + {\alpha_{2}} \right) \,
        \psi'({\alpha_{2}})\,\psi''({\alpha_{1}}) \right) } {8\,{\sigma_{12}}\,
    {\left( \psi'({\alpha_{1}}) + \psi'({\alpha_{2}}) -
        \left( {\alpha_{1}} + {\alpha_{2}} \right) \,\psi'({\alpha_{1}})
\,\psi'({\alpha_{2}})
        \right) }^2}  + \nonumber\\&&
        \frac{\psi'({\alpha_{1}})\,\left( -1 +
       \left( {\alpha_{1}} + {\alpha_{2}} \right) \,\psi'({\alpha_{1}}) \right)
\,
     \psi''({\alpha_{2}})} {8\,{\sigma_{12}}\,
    {\left( \psi'({\alpha_{1}}) + \psi'({\alpha_{2}}) -
        \left( {\alpha_{1}} + {\alpha_{2}} \right) \,\psi'({\alpha_{1}})
\,\psi'({\alpha_{2}})
        \right) }^2}\, ,\nonumber
\end{eqnarray*}
\begin{eqnarray}
 R_{33}&=&\frac{ \left( -2\,\psi'({\alpha_{1}})\,\psi'({\alpha_{2}}) +
       \left( {\alpha_{1}} + {\alpha_{2}} \right) \,\psi'({\alpha_{2}})\,
        \left( {\psi'({\alpha_{1}})}^2 - \psi''({\alpha_{1}}) \right)  +
      \psi''({\alpha_{1}}) \right) \,\psi''({\alpha_{2}})}{4\,
    {\left( \psi'({\alpha_{1}}) + \psi'({\alpha_{2}}) -
        \left( {\alpha_{1}} + {\alpha_{2}} \right) \,\psi'({\alpha_{1}})
\,\psi'({\alpha_{2}})
        \right) }^2}+ \nonumber\\&&
        \frac{{\psi'({\alpha_{1}})}^2\,{\psi'({\alpha_{2}})}^2}{4\,
    {\left( \psi'({\alpha_{1}}) + \psi'({\alpha_{2}}) -
        \left( {\alpha_{1}} + {\alpha_{2}} \right) \,\psi'({\alpha_{1}})
\,\psi'({\alpha_{2}})
        \right) }^2}\,.
        \end{eqnarray}
$\hfill \Box$
\end{proposition}

\begin{figure}
\begin{picture}(300,220)(0,0)
\put(-20,0){\includegraphics{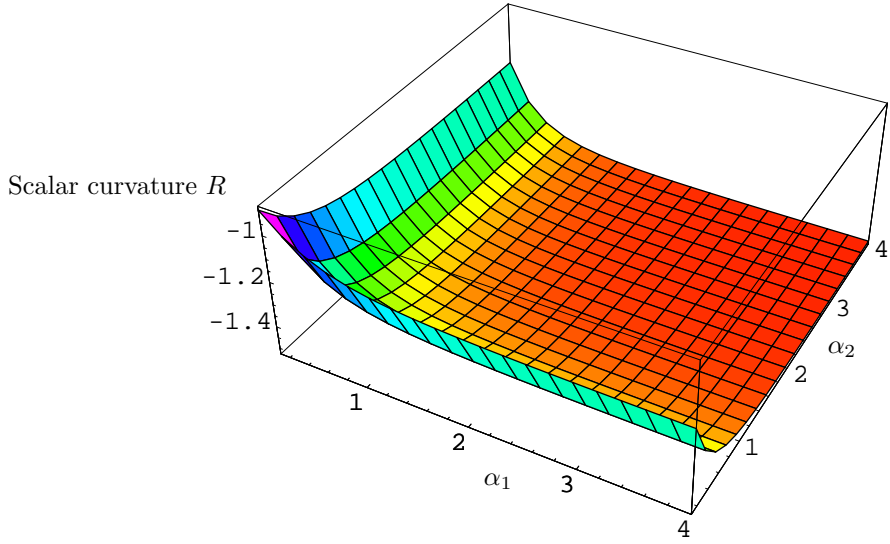}}
\put(310,80){$\alpha_2$}
\put(180,30){$\alpha_1$}
\put(0,140){Scalar curvature $R$}
\end{picture}
\caption{{\em The scalar curvature $R$ for the McKay bivariate gamma
3-manifold $M;$ the limiting value at the origin is $-\frac{1}{2}.$}}
\label{scalarmckay}
\end{figure}

\begin{proposition} The scalar curvature $R$ of $M$
 is given by:
\begin{eqnarray}
R&=& \frac{{\psi'(\alpha_{1})}^2\,\psi'(\alpha_{2}) +
     \psi'(\alpha_{1})\,{\psi'(\alpha_{2})}^2 +
     \alpha_{1}\,{\psi'(\alpha_{2})}^2\,\psi''(\alpha_{1})}{2\,
     {\left( \psi'(\alpha_{1}) + \psi'(\alpha_{2}) -
         \left( \alpha_{1} + \alpha_{2} \right) \,\psi'(\alpha_{1})\,
          \psi'(\alpha_{2}) \right) }^2} + \nonumber\\&&
  \frac{-\left( \psi'(\alpha_{2})\,\psi''(\alpha_{1}) \right)  -
     \psi'(\alpha_{1})\,\psi''(\alpha_{2})}{{\left( \psi'(\alpha_{1}) +
        \psi'(\alpha_{2}) - \left( \alpha_{1} + \alpha_{2} \right) \,
         \psi'(\alpha_{1})\,\psi'(\alpha_{2}) \right) }^2} +
         \nonumber\\&&
  \frac{ \alpha_{2}\,{\psi'(\alpha_{2})}^2\,\psi''(\alpha_{1}) +
     \alpha_{1}\,{\psi'(\alpha_{1})}^2\,\psi''(\alpha_{2}) +
    \alpha_{2}\,{\psi'(\alpha_{1})}^2\,\psi''(\alpha_{2})}{2\,
     {\left( \psi'(\alpha_{1}) + \psi'(\alpha_{2}) -
         \left( \alpha_{1} + \alpha_{2} \right) \,\psi'(\alpha_{1})\,
          \psi'(\alpha_{2}) \right) }^2} + \nonumber\\&&
  \frac{-\left( \alpha_{1}\,\psi''(\alpha_{1})\,\psi''(\alpha_{2}) \right)  -
     \alpha_{2}\,\psi''(\alpha_{1})\,\psi''(\alpha_{2})}{2\,
     {\left( \psi'(\alpha_{1}) + \psi'(\alpha_{2}) -
         \left( \alpha_{1} + \alpha_{2} \right) \,\psi'(\alpha_{1})\,
          \psi'(\alpha_{2}) \right) }^2}\, .
        \end{eqnarray}
This has limiting value $-\frac{1}{2}$ as $\alpha_{1}, \alpha_{2} \rightarrow 0.$
$\hfill \Box$
\end{proposition}
\begin{proposition} The sectional curvatures of $M$ are given by:
\begin{eqnarray}
 \varrho(1,2)&=& \frac{- \psi'({\alpha_{2}})\,
      \left( \psi'({\alpha_{1}}) +
        \left( {\alpha_{1}} + {\alpha_{2}} \right) \,\psi''({\alpha_{1}})
\right) }{4\,
    \left( -1 + \left( {\alpha_{1}} + {\alpha_{2}} \right) \,\psi'({\alpha_
{1}}) \right) \,
    \left( \psi'({\alpha_{1}}) + \psi'({\alpha_{2}}) -
      \left( {\alpha_{1}} + {\alpha_{2}} \right) \,\psi'({\alpha_{1}})
\,\psi'({\alpha_{2}})
      \right) } \,,\nonumber\\
 \varrho(1,3)&=&  \frac{-\left(
\psi'({\alpha_{1}})\,\psi'({\alpha_{2}}) +
      \left( \left( -3\,{\alpha_{1}} + {\alpha_{2}} \right) \,\psi'({\alpha_
{1}}) -
         4\,{{\alpha_{1}}}^2\,\psi''({\alpha_{1}}) \right) \,\psi''({\alpha_
{2}}) \right) }{4\,
    \left( \psi'({\alpha_{1}}) + \psi'({\alpha_{2}}) -
      \left( {\alpha_{1}} + {\alpha_{2}} \right) \,\psi'({\alpha_{1}})
\,\psi'({\alpha_{2}})
      \right) \,\left( -1 + \left( {\alpha_{2}} + {\alpha_{1}}\,
          \left( -3 + 4\,{\alpha_{1}}\,\psi'({\alpha_{1}}) \right)  \right)
\,\psi'({\alpha_{2}})
      \right) }\,,\nonumber\\
 \varrho(2,3)&=& \frac{- \psi'({\alpha_{1}})\,
      \left( \psi'({\alpha_{2}}) +
        \left( {\alpha_{1}} + {\alpha_{2}} \right) \,\psi''({\alpha_{2}})
\right) }{4\,
    \left( -1 + \left( {\alpha_{1}} + {\alpha_{2}} \right) \,\psi'({\alpha_
{2}}) \right) \,
    \left( \psi'({\alpha_{1}}) + \psi'({\alpha_{2}}) -
      \left( {\alpha_{1}} + {\alpha_{2}} \right) \,\psi'({\alpha_{1}})
\,\psi'({\alpha_{2}})
      \right) }\,.
      \end{eqnarray}
$\hfill \Box$
\end{proposition}

\begin{proposition}The mean curvatures $ \varrho(\lambda)\,(\lambda=1,2,3)$ are
given
by :
\begin{eqnarray}
\varrho(1)&=&\frac{-3\,{\alpha_{1}}\,{\psi'({\alpha_{1}})}^2\,\psi'({\alpha_
{2}})+
     {\alpha_{2}}\,{\psi'({\alpha_{1}})}^2\,\psi'({\alpha_{2}}) -
     3\,{\alpha_{1}}\,\psi'({\alpha_{1}})\,{\psi'({\alpha_{2}})}^2}{8\,
     \left( {\alpha_{2}} + {\alpha_{1}}\,\left( -3 + 4\,{\alpha_{1}}
\,\psi'({\alpha_{1}}) \right)  \right) \,
     {\left( \psi'({\alpha_{1}}) + \psi'({\alpha_{2}}) -
         \left( {\alpha_{1}} + {\alpha_{2}} \right) \,\psi'({\alpha_{1}})\,
          \psi'({\alpha_{2}}) \right) }^2} +\nonumber\\&&
  \frac{{\alpha_{2}}\,\psi'({\alpha_{1}})\,{\psi'({\alpha_{2}})}^2 +
     4\,{{\alpha_{1}}}^2\,{\psi'({\alpha_{1}})}^2\,{\psi'({\alpha_{2}})}^2 +
     3\,{\alpha_{1}}\,\psi'({\alpha_{2}})\,\psi''({\alpha_{1}})}{8\,
     \left( {\alpha_{2}} + {\alpha_{1}}\,\left( -3 + 4\,{\alpha_{1}}
\,\psi'({\alpha_{1}}) \right)  \right) \,
     {\left( \psi'({\alpha_{1}}) + \psi'({\alpha_{2}}) -
         \left( {\alpha_{1}} + {\alpha_{2}} \right) \,\psi'({\alpha_{1}})\,
          \psi'({\alpha_{2}}) \right) }^2} +\nonumber\\&&
  \frac{-\left( {\alpha_{2}}\,\psi'({\alpha_{2}})\,\psi''({\alpha_{1}})
\right)  -
     8\,{{\alpha_{1}}}^2\,\psi'({\alpha_{1}})\,\psi'({\alpha_{2}})\,
      \psi''({\alpha_{1}}) - 3\,{{\alpha_{1}}}^2\,{\psi'({\alpha_{2}})}^2\,
      \psi''({\alpha_{1}})}{8\,\left( {\alpha_{2}} +
       {\alpha_{1}}\,\left( -3 + 4\,{\alpha_{1}}\,\psi'({\alpha_{1}}) \right)
\right) \,
     {\left( \psi'({\alpha_{1}}) + \psi'({\alpha_{2}}) -
         \left( {\alpha_{1}} + {\alpha_{2}} \right) \,\psi'({\alpha_{1}})\,
          \psi'({\alpha_{2}}) \right) }^2} +\nonumber\\&&
  \frac{-2\,{\alpha_{1}}\,{\alpha_{2}}\,{\psi'({\alpha_{2}})}^2\,\psi''({\alpha_
{1}}) +
     {{\alpha_{2}}}^2\,{\psi'({\alpha_{2}})}^2\,\psi''({\alpha_{1}}) +
     4\,{{\alpha_{1}}}^3\,\psi'({\alpha_{1}})\,{\psi'({\alpha_{2}})}^2\,
      \psi''({\alpha_{1}})}{8\,\left( {\alpha_{2}} +
       {\alpha_{1}}\,\left( -3 + 4\,{\alpha_{1}}\,\psi'({\alpha_{1}}) \right)
\right) \,
     {\left( \psi'({\alpha_{1}}) + \psi'({\alpha_{2}}) -
         \left( {\alpha_{1}} + {\alpha_{2}} \right) \,\psi'({\alpha_{1}})\,
          \psi'({\alpha_{2}}) \right) }^2} +\nonumber\\&&
  \frac{4\,{{\alpha_{1}}}^2\,{\alpha_{2}}\,\psi'({\alpha_{1}})\,{\psi'({\alpha_
{2}})}^2\,
      \psi''({\alpha_{1}}) + 3\,{\alpha_{1}}\,\psi'({\alpha_{1}})
\,\psi''({\alpha_{2}}) -
     {\alpha_{2}}\,\psi'({\alpha_{1}})\,\psi''({\alpha_{2}})}{8\,
     \left( {\alpha_{2}} + {\alpha_{1}}\,\left( -3 + 4\,{\alpha_{1}}
\,\psi'({\alpha_{1}}) \right)  \right) \,
     {\left( \psi'({\alpha_{1}}) + \psi'({\alpha_{2}}) -
         \left( {\alpha_{1}} + {\alpha_{2}} \right) \,\psi'({\alpha_{1}})\,
          \psi'({\alpha_{2}}) \right) }^2} +\nonumber\\&&
  \frac{-3\,{{\alpha_{1}}}^2\,{\psi'({\alpha_{1}})}^2\,\psi''({\alpha_{2}}) -
     2\,{\alpha_{1}}\,{\alpha_{2}}\,{\psi'({\alpha_{1}})}^2\,\psi''({\alpha_
{2}}) +
     {{\alpha_{2}}}^2\,{\psi'({\alpha_{1}})}^2\,\psi''({\alpha_{2}})}{8\,
     \left( {\alpha_{2}} + {\alpha_{1}}\,\left( -3 + 4\,{\alpha_{1}}
\,\psi'({\alpha_{1}}) \right)  \right) \,
     {\left( \psi'({\alpha_{1}}) + \psi'({\alpha_{2}}) -
         \left( {\alpha_{1}} + {\alpha_{2}} \right) \,\psi'({\alpha_{1}})\,
          \psi'({\alpha_{2}}) \right) }^2} +\nonumber\\&&
  \frac{{{\alpha_{1}}}^2\,\psi''({\alpha_{1}})\,\psi''({\alpha_{2}}) -
     {{\alpha_{1}}}^3\,\psi'({\alpha_{1}})\,\psi''({\alpha_{1}})
\,\psi''({\alpha_{2}}) -
     {{\alpha_{1}}}^2\,{\alpha_{2}}\,\psi'({\alpha_{1}})\,\psi''({\alpha_{1}})\,
      \psi''({\alpha_{2}})}{2\,\left( {\alpha_{2}} +
       {\alpha_{1}}\,\left( -3 + 4\,{\alpha_{1}}\,\psi'({\alpha_{1}}) \right)
\right) \,
     {\left( \psi'({\alpha_{1}}) + \psi'({\alpha_{2}}) -
         \left( {\alpha_{1}} + {\alpha_{2}} \right) \,\psi'({\alpha_{1}})\,
          \psi'({\alpha_{2}}) \right) }^2} \,,\nonumber\\
 \varrho(2)&=& \frac{\psi'({\alpha_{2}})\,\left(
\psi'({\alpha_{1}})\,
        \left( \psi'({\alpha_{1}}) + \psi'({\alpha_{2}}) \right)  -
       \psi''({\alpha_{1}}) + \left( {\alpha_{1}} + {\alpha_{2}} \right) \,
        \psi'({\alpha_{2}})\,\psi''({\alpha_{1}}) \right) }
     {8\,{\left( \psi'({\alpha_{1}}) + \psi'({\alpha_{2}}) -
        \left( {\alpha_{1}} + {\alpha_{2}} \right) \,\psi'({\alpha_{1}})
\,\psi'({\alpha_{2}})
        \right) }^2} + \nonumber\\&&
         \frac{\psi'({\alpha_{1}})\,\left( -1 +
       \left( {\alpha_{1}} + {\alpha_{2}} \right) \,\psi'({\alpha_{1}}) \right)
\,
     \psi''({\alpha_{2}})}
     {8\,{\left( \psi'({\alpha_{1}}) + \psi'({\alpha_{2}}) -
        \left( {\alpha_{1}} + {\alpha_{2}} \right) \,\psi'({\alpha_{1}})
\,\psi'({\alpha_{2}})
        \right) }^2}\, , \nonumber\\
 \varrho(3)&=&\frac{ \left( -2\,\psi'({\alpha_{1}})\,\psi'({\alpha_{2}}) +
       \left( {\alpha_{1}} + {\alpha_{2}} \right) \,\psi'({\alpha_{2}})\,
        \left( {\psi'({\alpha_{1}})}^2 - \psi''({\alpha_{1}}) \right)  +
       \psi''({\alpha_{1}}) \right) \,\psi''({\alpha_{2}})}
       {8\,\psi'({\alpha_{2}})\,
    {\left( \psi'({\alpha_{1}}) + \psi'({\alpha_{2}}) -
        \left( {\alpha_{1}} + {\alpha_{2}} \right) \,\psi'({\alpha_{1}})
\,\psi'({\alpha_{2}})
        \right) }^2}+ \nonumber\\&&
 \frac{{\psi'({\alpha_{1}})}^2\,{\psi'({\alpha_{2}})}^2 }
       {8\,\psi'({\alpha_{2}})\,
    {\left( \psi'({\alpha_{1}}) + \psi'({\alpha_{2}}) -
        \left( {\alpha_{1}} + {\alpha_{2}} \right) \,\psi'({\alpha_{1}})
\,\psi'({\alpha_{2}})
        \right) }^2}\, .
        \end{eqnarray}
$\hfill \Box$
         \end{proposition}

\section{Submanifolds of the Mckay 3-manifold $M$}
We consider three submanifolds
$M_{1}, M_{2}$ and $ M_{3}$ of the 3-manifold $M$ of Mckay bivariate
gamma distributions (\ref{mckaymodel}) $
f(x,y;\alpha_{1},\sigma_{12},\alpha_{2})$, where we use the
coordinate system $ (\alpha_{1},\sigma_{12},\alpha_{2}).$
These submanifolds have dimension 2 and so it follows that
the scalar curvature is twice the Gaussian curvature,
$R=2K.$ Recall from above that the correlation is given by
$$\rho = \sqrt{\frac{\alpha_{1}}{\alpha_{1}+\alpha_{2}}}.$$
In the cases of $M_{1}$ and $ M_{2}$ the scalar curvature can be
shown as a function only of $\rho$.

\subsection{Submanifold $ M_{1}\subset M$: $\alpha_{1}=1$}
The distributions are of form:
 \begin{eqnarray}
f(x,y;1,\sigma_{12},\alpha_{2}) =
\frac{(\frac{1}{\sigma_{12}})^{\frac{1+\alpha_{2}}{2}}(y-x)^{\alpha_{2}-1}
e^{-\sqrt{\frac{1}{\sigma_{12}}}y}}{\Gamma(\alpha_{2})} \ ,
\end{eqnarray}
defined on $ 0<x<y<\infty $ with parameters $
\sigma_{12},\alpha_{2}>0.$  The correlation coefficient and
marginal functions, of $X$ and $Y$ are given by :
\begin{eqnarray}
 \rho(X,Y)&=&\frac{1}{\sqrt{1+\alpha_{2}}}  \\
 f_{X}(x) &= & \frac{1}{\sqrt{\sigma_{12}}}\, e^{-\frac{1}{\sqrt{\sigma_{12}}}
x},\quad x
>0 \\
 f_{Y}(y)&=&\frac{(\frac{1}{\sqrt{\sigma_{12}}})^{(1+\alpha_{2})}
y^{\alpha_{2}}e^ {-\frac{1}{\sqrt{\sigma_{12}}}
y}}{\alpha_{2}\,\Gamma(\alpha_{2})}, \quad y>0
\end{eqnarray}
So here we have $ \alpha_{2}=\frac{1-\rho^2}{\rho^2},$ which in practice
would give a measure of the variability not due to the correlation.

\begin{proposition}
The metric tensor $[g_{ij}]$ and
its inverse $ [g^{ij}]$ are as follows :
\begin{eqnarray}
 G=[g_{ij}]=\left[ \begin{array}{ccc}
 \frac{1 + {\alpha_{2}}}{4\,{{\sigma_{12}}}^2} & \frac{1}{2\,{\sigma_{12}}} \\
 \frac{1}{2\,{\sigma_{12}}}             & \psi'(\alpha_{2})
\end{array} \right]
\label{1mckaymetric}
\end{eqnarray}
\begin{eqnarray}
 G^{-1}=[g^{ij}]=\left[ \begin{array}{ccc}
 \frac{4\,{{\sigma_{12}}}^2\,\psi'(\alpha_{2})}{-1 + \left( 1 + {\alpha_{2}}
\right)\,\psi'(\alpha_{2})}
  & \frac{-2\,{\sigma_{12}}}{-1 + \left( 1 + {\alpha_{2}} \right)\,\psi'(\alpha_
{2})}\\
 \frac{-2\,{\sigma_{12}}}{-1 + \left( 1 + {\alpha_{2}} \right)\,\psi'(\alpha_
{2})}
 & \frac{1 + {\alpha_{2}}}{-1 + \left( 1 + {\alpha_{2}} \right) \,\psi'(\alpha_
{2})}
\end{array} \right]
\end{eqnarray}
$\hfill \Box$
         \end{proposition}

\begin{proposition}The Christoffel symbols of $M_{1}$ are
\begin{eqnarray}
 \Gamma^{1}_{11}&=& \frac{-4 + \frac{1}
     {-1 + \left( 1 + {\alpha_{2}} \right) \,\psi'(\alpha_{2})}}
    {4\,{\sigma_{12}}} \,,\nonumber\\
 \Gamma^{1}_{12}&=& \Gamma^{1}_{21}= \frac{\psi'(\alpha_{2})}
  {-2 + 2\,\left( 1 + {\alpha_{2}} \right) \,\psi'(\alpha_{2})}\,,
  \nonumber\\
 \Gamma^{1}_{22}&=& -\left(
\frac{{\sigma_{12}}\,\psi''(\alpha_{2})}
    {-1 + \left( 1 + {\alpha_{2}} \right) \,\psi'(\alpha_{2})}
    \right)\, ,\nonumber\\
 \Gamma^{2}_{11}&=& \frac{-\left( 1 + {\alpha_{2}} \right) }
  {8\,{{\sigma_{12}}}^2\,\left( -1 + \left( 1 + {\alpha_{2}} \right) \,
       \psi'(\alpha_{2}) \right) }\, ,\nonumber\\
 \Gamma^{2}_{12}&=& \Gamma^{2}_{21}=
\frac{-1}{4\,{\sigma_{12}}\,\left( -1 +
      \left( 1 + {\alpha_{2}} \right) \,\psi'(\alpha_{2}) \right) }\,,
      \nonumber\\
 \Gamma^{2}_{22}&=& \frac{\left( 1 + {\alpha_{2}} \right)
\,\psi''(\alpha_{2})}
  {-2 + 2\,\left( 1 + {\alpha_{2}} \right) \,\psi'(\alpha_{2})} \,.
  \end{eqnarray}
$\hfill \Box$
         \end{proposition}

\begin{proposition}
The curvature tensor of $M_{1}$ is given by
\begin{eqnarray}
 R_{1212}&=& \frac{-\left( \psi'({\alpha_{2}}) + \left( 1 + {\alpha_{2}}
\right) \,\psi''({\alpha_{2}}) \right) }
  {16\,{{\sigma_{12}}}^2\,\left( -1 + \left( 1 + {\alpha_{2}} \right)
\,\psi'({\alpha_{2}}) \right) } \, ,
       \end{eqnarray}
       while the other independent components are zero.

  By contraction we obtain:

  Ricci tensor:
  \begin{eqnarray}
      R_{11}&=& \frac{\left( 1 + {\alpha_{2}} \right) \,\left( \psi'({\alpha_
{2}}) +
      \left( 1 + {\alpha_{2}} \right) \,\psi''({\alpha_{2}}) \right) }{16\,
{{\sigma_{12}}}^2\,
    {\left( -1 + \left( 1 + {\alpha_{2}} \right) \,\psi'({\alpha_{2}}) \right) }
^2}\,,\nonumber\\
       R_{12}&=&\frac{\psi'({\alpha_{2}}) + \left( 1 + {\alpha_{2}} \right)
\,\psi''({\alpha_{2}})}
  {8\,{\sigma_{12}}\,{\left( -1 + \left( 1 + {\alpha_{2}} \right)
\,\psi'({\alpha_{2}}) \right) }^2} \,,\nonumber\\
      R_{22}&=& \frac{\psi'({\alpha_{2}})\,\left( \psi'({\alpha_{2}}) +
      \left( 1 + {\alpha_{2}} \right) \,\psi''({\alpha_{2}}) \right) }{4\,
    {\left( -1 + \left( 1 + {\alpha_{2}} \right) \,\psi'({\alpha_{2}}) \right) }
^2}\,.
         \end{eqnarray}

\begin{figure}
\begin{picture}(300,220)(0,0)
\put(-20,20){\includegraphics{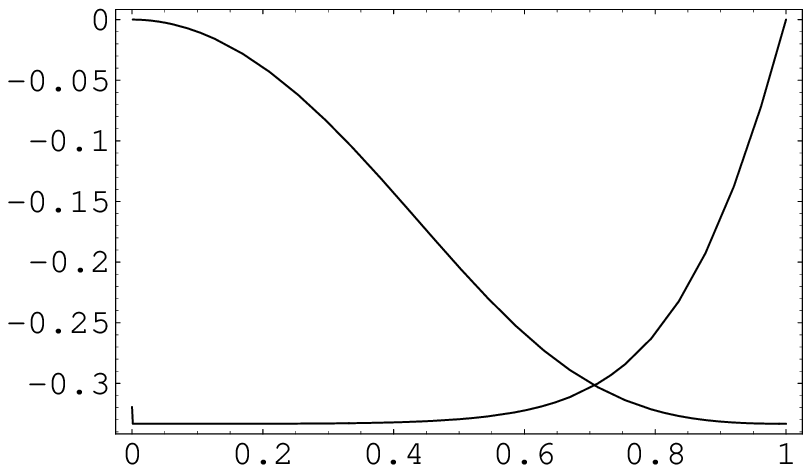}}
\put(284,110){$M_1$}
\put(150,140){$M_2$}
\put(40,180){Scalar
curvature $R$} \put(180,20){Correlation $\rho$}
\end{picture}
\caption{{\em The scalar curvature $R$ as a function of
correlation $\rho$ for McKay submanifolds: $M_1$ ($M$ with
$\alpha_1=1$) where $R$ increases from $-\frac{1}{3}$ to $0,$ and
$M_2$ ($M$ with $\alpha_2=1$) where $R$ decreases from $0$ to
$-\frac{1}{3}.$}} \label{RrhoM1M3}
\end{figure}

Scalar curvature:
\begin{eqnarray}
 R&=& \frac{\psi'({\alpha_{2}}) + \left( 1 + {\alpha_{2}} \right)
\,\psi''({\alpha_{2}})}
  {2\,{\left( -1 + \left( 1 + {\alpha_{2}} \right) \,\psi'({\alpha_{2}})
\right) }^2}\,.
\end{eqnarray}

$\hfill \Box$
         \end{proposition}

\subsection{Submanifold $ M_{2}\subset M$: $\alpha_{2}=1$}
 The distributions are of form:
 \begin{eqnarray}
f(x,y;\alpha_{1},\sigma_{12},1) =
\frac{(\frac{\alpha_{1}}{\sigma_{12}})^{\frac{\alpha_{1}+1}{2}}x^{\alpha_{1}-1}
e^{-\sqrt{\frac{\alpha_{1}}{\sigma_{12}}}y}}{\Gamma(\alpha_{1})} \
,
\end{eqnarray}
defined on $ 0<x<y<\infty $ with parameters $
\alpha_1,\sigma_{12}>0$. The correlation coefficient and marginal
functions, of $X$ and $Y$ are given by :
\begin{eqnarray}
 \rho(X,Y)&=& \sqrt{\frac{\alpha_{1}}{1+\alpha_{1}}}\\
 f_{X}(x) &=&  \frac{(\frac{\alpha_{1}}{\sigma_{12}})^{\frac{\alpha_{1}}{2}} x^
{\alpha_{1}-1}e^{-\sqrt{\frac{\alpha_{1}}{\sigma_{12}}}
x}}{\Gamma(\alpha_{1})},\quad x
>0 \\
 f_{Y}(y)&=&\frac{(\frac{\alpha_{1}}{\sigma_{12}})^{\frac{(\alpha_{1}+1)}{2}}
y^{\alpha_{1}}e^ {-\sqrt{\frac{\alpha_{1}}{\sigma_{12}}}
y}}{\alpha_{1}\,\Gamma(\alpha_{1})}, \quad y>0
\end{eqnarray}
Here we have $ \alpha_{1}=\frac{\rho^2}{1-\rho^2}.$

\begin{proposition}
The metric tensor $[g_{ij}]$ and its inverse $ [g^{ij}]$ are as
follows :
\begin{eqnarray}
 G=[g_{ij}]=\left[ \begin{array}{ccc}
\frac{1 - 3\,{\alpha_{1}}}{4\,{{\alpha_{1}}}^2} +
\psi'(\alpha_{1})
& \frac{-1 + {\alpha_{1}}}{4\,{\alpha_{1}}\,{\sigma_{12}}} \\
 \frac{-1 + {\alpha_{1}}}{4\,{\alpha_{1}}\,{\sigma_{12}}} & \frac{1 + {\alpha_
{1}}}{4\,{{\sigma_{12}}}^2}
\end{array} \right]
\label{3mckaymetric}
\end{eqnarray}
\begin{eqnarray}
 G^{-1}=[g^{ij}]=\left[ \begin{array}{ccc}
 \frac{1 + {\alpha_{1}}}{-1 + \left( 1 + {\alpha_{1}} \right)
 \,\psi'(\alpha_{1})} &
   -\left( \frac{\left( -1 + {\alpha_{1}} \right) \,{\sigma_{12}}}
      {{\alpha_{1}}\,\left( -1 + \left( 1 + {\alpha_{1}} \right) \,\psi'(\alpha_
{1}) \right) }
      \right)  \\  -\left( \frac{\left( -1 + {\alpha_{1}} \right) \,{\sigma_
{12}}}
      {{\alpha_{1}}\,\left( -1 + \left( 1 + {\alpha_{1}} \right) \,\psi'(\alpha_
{1}) \right) }
      \right) & \frac{{{\sigma_{12}}}^2\,\left( 1 +
        {\alpha_{1}}\,\left( -3 + 4\,{\alpha_{1}}\,\psi'(\alpha_{1}) \right)
\right) }{{{\alpha_{1}}}^2\,
      \left( -1 + \left( 1 + {\alpha_{1}} \right) \,\psi'(\alpha_{1}) \right) }
\end{array} \right]
\end{eqnarray}
$\hfill \Box$
         \end{proposition}

\begin{proposition}
The Christoffel symbols are
\begin{eqnarray}
 \Gamma^{1}_{11}&=&\frac{-1 + {\alpha_{1}}\,\left( 3 +
4\,{\alpha_{1}}\,\left( 1 + {\alpha_{1}} \right) \,
        \psi''(\alpha_{1}) \right) }{8\,{{\alpha_{1}}}^2\,
    \left( -1 + \left( 1 + {\alpha_{1}} \right) \,\psi'(\alpha_{1}) \right) }
   \, ,\nonumber\\
 \Gamma^{1}_{12}&=& \Gamma^{1}_{21}=\frac{-\left( -1 + {\alpha_{1}}
\right) }
  {8\,{\alpha_{1}}\,{\sigma_{12}}\,\left( -1 + \left( 1 + {\alpha_{1}} \right)
\,\psi'(\alpha_{1}) \right)
      } \,,\nonumber\\
\Gamma^{1}_{22}&=&\frac{-\left( 1 + {\alpha_{1}} \right) }
  {8\,{{\sigma_{12}}}^2\,\left( -1 + \left( 1 + {\alpha_{1}} \right)
\,\psi'(\alpha_{1}) \right) }
 \, ,\nonumber\\
 \Gamma^{2}_{11}&=& \frac{{\sigma_{12}}\,\left( -1 +
{\alpha_{1}}\,\left( -3 + 8\,\psi'(\alpha_{1}) -
         4\,\left( -1 + {\alpha_{1}} \right) \,{\alpha_{1}}\,\psi''(\alpha_{1})
\right)  \right) }
    {8\,{{\alpha_{1}}}^3\,\left( -1 + \left( 1 + {\alpha_{1}} \right)
\,\psi'(\alpha_{1}) \right) }
   \, ,\nonumber\\
 \Gamma^{2}_{12}&=& \frac{1 + {\alpha_{1}}\,\left( -3 +
4\,{\alpha_{1}}\,\psi'(\alpha_{1}) \right) }
  {8\,{{\alpha_{1}}}^2\,\left( -1 + \left( 1 + {\alpha_{1}} \right)
\,\psi'(\alpha_{1}) \right) }
  ,\nonumber\\
 \Gamma^{2}_{22}&=&\frac{-8 + \frac{-1 + {\alpha_{1}}}
     {{\alpha_{1}}\,\left( -1 + \left( 1 + {\alpha_{1}} \right) \,\psi'(\alpha_
{1}) \right) }}
    {8\,{\sigma_{12}}} \,.
    \end{eqnarray}
$\hfill \Box$
         \end{proposition}

\begin{proposition}
 The curvature tensor is given by
\begin{eqnarray}
 R_{1212}&=&  \frac{-\left( \psi'({\alpha_{1}}) + \left( 1 + {\alpha_{1}}
\right) \,\psi''({\alpha_{1}}) \right) }
  {16\,{{\sigma_{12}}}^2\,\left( -1 + \left( 1 + {\alpha_{1}} \right)
\,\psi'({\alpha_{1}}) \right) } .
\end{eqnarray}
while the other independent components are zero.

 By contraction we obtain:\\
 Ricci tensor:
 \begin{eqnarray}
 R_{11}&=& \frac{\left( 1 + {\alpha_{1}}\,\left( -3 + 4\,{\alpha_{1}}
\,\psi'({\alpha_{1}}) \right)  \right) \,
    \left( \psi'({\alpha_{1}}) + \left( 1 + {\alpha_{1}} \right)
\,\psi''({\alpha_{1}}) \right) }
    {16\,{{\alpha_{1}}}^2\,{\left( -1 + \left( 1 + {\alpha_{1}} \right)
\,\psi'({\alpha_{1}}) \right) }^2}
    \,,\nonumber\\
     R_{21}&=& \frac{\left( -1 + {\alpha_{1}} \right) \,\left( \psi'({\alpha_
{1}}) +
      \left( 1 + {\alpha_{1}} \right) \,\psi''({\alpha_{1}}) \right) }{16\,
{\alpha_{1}}\,{\sigma_{12}}\,
    {\left( -1 + \left( 1 + {\alpha_{1}} \right) \,\psi'({\alpha_{1}}) \right) }
^2}
    \,,\nonumber\\
     R_{22}&=& \frac{\left( 1 + {\alpha_{1}} \right) \,\left( \psi'({\alpha_
{1}}) +
      \left( 1 + {\alpha_{1}} \right) \,\psi''({\alpha_{1}}) \right) }{16\,
{{\sigma_{12}}}^2\,
    {\left( -1 + \left( 1 + {\alpha_{1}} \right) \,\psi'({\alpha_{1}}) \right) }
^2} \,.
    \end{eqnarray}
Scalar curvature:
\begin{eqnarray}
 R&=& \frac{\psi'({\alpha_{1}}) + \left( 1 + {\alpha_{1}} \right)
\,\psi''({\alpha_{1}})}
  {2\,{\left( -1 + \left( 1 + {\alpha_{1}} \right) \,\psi'({\alpha_{1}})
\right) }^2}\,.
    \end{eqnarray}

$\hfill \Box$
         \end{proposition}

  \subsection{Submanifold $ M_{3}\subset M$: $\sigma_{12}=1$}
Here the distributions have unit covariance and are of form:
 \begin{eqnarray}
f(x,y;\alpha_{1},1,\alpha_{2}) =
\frac{(\alpha_{1})^{\frac{(\alpha_{1}+\alpha_{2})}{2}}x^{\alpha_{1}-1}(y-x)^
{\alpha_{2}-1}
e^{-\sqrt{\alpha_{1}}y}}{\Gamma(\alpha_{1})\Gamma(\alpha_{2})} \ ,
\end{eqnarray}
defined on $ 0<x<y<\infty $ with parameters $
\alpha_{1},\alpha_{2}>0$.  The correlation coefficient and
marginal functions, of $X$ and $Y$ are given by :
\begin{eqnarray}
 \rho(X,Y)&=&\sqrt{\frac{\alpha_{1}}{\alpha_{1}+\alpha_{2}}}\\
 f_{X}(x) &= & \frac{{\sqrt{\alpha_{1}}}^{{\alpha_{1}}} x^
{\alpha_{1}-1}e^{-\sqrt{\alpha_{1}} x}}{\Gamma(\alpha_{1})},\quad
x
>0 \\
 f_{Y}(y)&=&\frac{{\sqrt{\alpha_{1}}}^{(\alpha_{1}+\alpha_{2})}
y^{(\alpha_{1}+\alpha_{2})-1}e^ {-\sqrt{\alpha_{1}}
y}}{\Gamma(\alpha_{1}+\alpha_{2})}, \quad y>0
\end{eqnarray}

\begin{proposition}
The metric tensor $[g_{ij}]$ and
its inverse $ [g^{ij}]$ are as follows:
\begin{eqnarray}
 G=[g_{ij}]=\left[ \begin{array}{ccc}
 \frac{-3\,{\alpha_{1}} + {\alpha_{2}}}{4\,{{\alpha_{1}}}^2} +
\psi'(\alpha_{1}) & \frac{-1}{2\,{\alpha_{1}}} \\
\frac{-1}{2\,{\alpha_{1}}} & \psi'(\alpha_{2})
\end{array} \right]
\label{2mckaymetric}
 \end{eqnarray}
\begin{eqnarray}
G^{-1}=[g^{ij}]=\left[ \begin{array}{ccc}
 \frac{4\,{{\alpha_{1}}}^2\,\psi'(\alpha_{2})}
    {-1 + \left( {\alpha_{2}} + {\alpha_{1}}\,\left( -3 + 4\,{\alpha_{1}}
\,\psi'(\alpha_{1}) \right)
         \right) \,\psi'(\alpha_{2})} &
   \frac{2\,{\alpha_{1}}}{-1 + \left( {\alpha_{2}} + {\alpha_{1}}\,
          \left( -3 + 4\,{\alpha_{1}}\,\psi'(\alpha_{1}) \right)  \right) \,
       \psi'(\alpha_{2})}\\
  \frac{2\,{\alpha_{1}}}
    {-1 + \left( {\alpha_{2}} + {\alpha_{1}}\,\left( -3 + 4\,{\alpha_{1}}
\,\psi'(\alpha_{1}) \right)
         \right) \,\psi'(\alpha_{2})} &
   \frac{{\alpha_{2}} + {\alpha_{1}}\,\left( -3 + 4\,{\alpha_{1}}\,\psi'(\alpha_
{1}) \right) }
    {-1 + \left( {\alpha_{2}} + {\alpha_{1}}\,\left( -3 + 4\,{\alpha_{1}}
\,\psi'(\alpha_{1}) \right)
         \right) \,\psi'(\alpha_{2})}
\end{array} \right]
\end{eqnarray}
$\hfill \Box$
         \end{proposition}

\begin{proposition}The Christoffel symbols are
\begin{eqnarray}
 \Gamma^{1}_{11}&=& \frac{3 + 2\,\psi'(\alpha_{2})\,
     \left( 3\,{\alpha_{1}} - 2\,{\alpha_{2}} + 4\,{{\alpha_{1}}}^3
\,\psi''(\alpha_{1}) \right) }{4\,{\alpha_{1}}\,
    \left( -1 + \left( {\alpha_{2}} + {\alpha_{1}}\,\left( -3 + 4\,{\alpha_{1}}
\,\psi'(\alpha_{1}) \right)
             \right) \,\psi'(\alpha_{2}) \right) } \,,\nonumber\\
 \Gamma^{1}_{22}&=& \frac{{\alpha_{1}}\,\psi''(\alpha_{2})}
  {-1 + \left( {\alpha_{2}} + {\alpha_{1}}\,\left( -3 + 4\,{\alpha_{1}}
\,\psi'(\alpha_{1}) \right)  \right) \,
     \psi'(\alpha_{2})}\, ,\nonumber\\
  \Gamma^{1}_{12}&=& \Gamma^{1}_{21}=\frac{\psi'(\alpha_{2})}
  {-2 + 2\,\left( {\alpha_{2}} + {\alpha_{1}}\,\left( -3 + 4\,{\alpha_{1}}
\,\psi'(\alpha_{1}) \right)  \right)
       \,\psi'(\alpha_{2})} \,, \nonumber\\
 \Gamma^{2}_{11}&=& \frac{-{\alpha_{2}} + {\alpha_{1}}\,\left( -3 +
12\,{\alpha_{1}}\,\psi'(\alpha_{1}) +
       8\,{{\alpha_{1}}}^2\,\psi''(\alpha_{1}) \right) }{8\,{{\alpha_{1}}}^2\,
    \left( -1 + \left( {\alpha_{2}} + {\alpha_{1}}\,\left( -3 + 4\,{\alpha_{1}}
\,\psi'(\alpha_{1}) \right)
             \right) \,\psi'(\alpha_{2}) \right) }\, ,\nonumber\\
  \Gamma^{2}_{22}&=&\frac{\left( {\alpha_{2}} + {\alpha_{1}}\,\left( -3 +
4\,{\alpha_{1}}\,\psi'(\alpha_{1}) \right)  \right) \,
    \psi''(\alpha_{2})}{-2 +
    2\,\left( {\alpha_{2}} + {\alpha_{1}}\,\left( -3 + 4\,{\alpha_{1}}
\,\psi'(\alpha_{1}) \right)  \right) \,
     \psi'(\alpha_{2})}\,, \nonumber \\
  \Gamma^{2}_{12}&=& \Gamma^{2}_{21}= \frac{1}{4\,{\alpha_{1}}\,\left( -1 +
\left( {\alpha_{2}} +
         {\alpha_{1}}\,\left( -3 + 4\,{\alpha_{1}}\,\psi'(\alpha_{1}) \right)
\right) \,
       \psi'(\alpha_{2}) \right) }\, .
       \end{eqnarray}
$\hfill \Box$
         \end{proposition}

\begin{proposition}
The curvature tensor is given by
\begin{eqnarray}
 R_{1212}&=& \frac{-\psi'({\alpha_{2}}) + \left( -{\alpha_{2}}+
       {\alpha_{1}}\,\left( -3 + 12\,{\alpha_{1}}\,\psi'({\alpha_{1}}) +
          8\,{{\alpha_{1}}}^2\,\psi''({\alpha_{1}}) \right)  \right)
\,\psi''({\alpha_{2}})}{16\,
    {{\alpha_{1}}}^2\,\left( -1 + \left( {\alpha_{2}} + {\alpha_{1}}\,\left( -3
+ 4\,{\alpha_{1}}\,\psi'({\alpha_{1}}) \right)
         \right) \,\psi'({\alpha_{2}}) \right) }\, ,
         \end{eqnarray}
       while the other independent components are zero.

\begin{figure}
\begin{picture}(300,220)(0,0)
\put(-10,20){\includegraphics{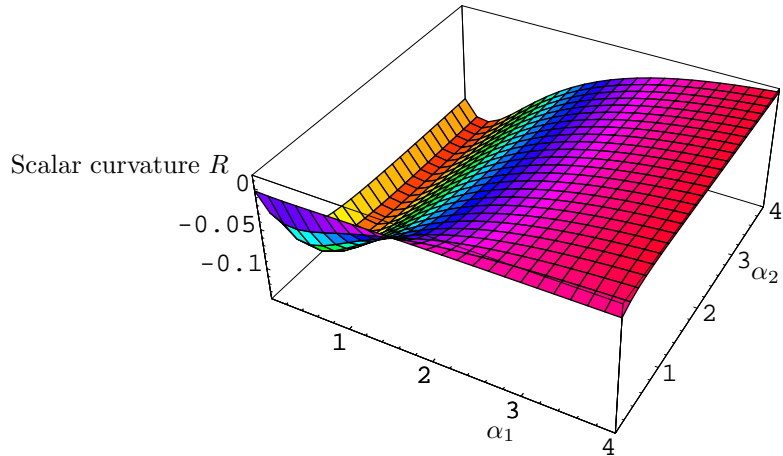}}
\put(300,100){$\alpha_2$}
\put(200,40){$\alpha_1$}
\put(20,140){Scalar curvature $R$}
\end{picture}
\caption{{\em The scalar curvature $R$ for McKay submanifold
$M_3$, ($M$ with $\sigma_{12}=1$).}} \label{RalphM2}
\end{figure}

  By contraction we obtain:\\
  Ricci tensor:
   \begin{eqnarray}
         R_{11}&=& \frac{ \left( {\alpha_{2}} + {\alpha_{1}}\,\left( -3 + 4\,
{\alpha_{1}}\,
         \psi'({\alpha_{1}}) \right)  \right) \,
      \left( \left( -{\alpha_{2}} + {\alpha_{1}}\,\left( -3 + 12\,{\alpha_{1}}
\,\psi'({\alpha_{1}}) +
              8\,{{\alpha_{1}}}^2\,\psi''({\alpha_{1}}) \right)  \right)
\,\psi''({\alpha_{2}}) \right)  }
        {- 16\,{{\alpha_{1}}}^2\,{\left( -1 +
        \left( {\alpha_{2}} + {\alpha_{1}}\,\left( -3 + 4\,{\alpha_{1}}
\,\psi'({\alpha_{1}}) \right)  \right) \,
         \psi'({\alpha_{2}}) \right) }^2} \nonumber\\&&
         + \frac{ \psi'({\alpha_{2}}) \,\left( {\alpha_{2}} + {\alpha_{1}}
\,\left( -3 + 4\,{\alpha_{1}}\,
         \psi'({\alpha_{1}}) \right)  \right) \,
      \left( \psi'({\alpha_{2}})\right)  }
        {16\,{{\alpha_{1}}}^2\,{\left( -1 +
        \left( {\alpha_{2}} + {\alpha_{1}}\,\left( -3 + 4\,{\alpha_{1}}
\,\psi'({\alpha_{1}}) \right)  \right) \,
         \psi'({\alpha_{2}}) \right) }^2}\, ,\nonumber\\
     R_{12}&=& \frac{-\psi'({\alpha_{2}}) + \left( -{\alpha_{2}} +
       {\alpha_{1}}\,\left( -3 + 12\,{\alpha_{1}}\,\psi'({\alpha_{1}}) +
          8\,{{\alpha_{1}}}^2\,\psi''({\alpha_{1}}) \right)  \right)
\,\psi''({\alpha_{2}})}{8\,
    {\alpha_{1}}\,{\left( -1 + \left( {\alpha_{2}} + {\alpha_{1}}\,\left( -3 + 4
\,{\alpha_{1}}\,\psi'({\alpha_{1}}) \right)  \right)
           \,\psi'({\alpha_{2}}) \right) }^2}\, ,\nonumber\\
      R_{22}&=&\frac{\psi'({\alpha_{2}})\,\left( \psi'({\alpha_{2}}) +
      \left( {\alpha_{2}} + {\alpha_{1}}\,\left( 3 -
            4\,{\alpha_{1}}\,\left( 3\,\psi'({\alpha_{1}}) + 2\,{\alpha_{1}}
\,\psi''({\alpha_{1}}) \right)
            \right)  \right) \,\psi''({\alpha_{2}}) \right) }{4\,
    {\left( -1 + \left( {\alpha_{2}} + {\alpha_{1}}\,\left( -3 + 4\,{\alpha_{1}}
\,\psi'({\alpha_{1}}) \right)  \right) \,
         \psi'({\alpha_{2}}) \right) }^2}\, .
         \end{eqnarray}
 Scalar curvature:
\begin{eqnarray}
 R&=& \frac{\psi'({\alpha_{2}}) + \left( {\alpha_{2}} +
       {\alpha_{1}}\,\left( 3 - 4\,{\alpha_{1}}\,\left( 3\,\psi'({\alpha_{1}}) +
             2\,{\alpha_{1}}\,\psi''({\alpha_{1}}) \right)  \right)  \right)
\,\psi''({\alpha_{2}})}
    {2\,{\left( -1 + \left( {\alpha_{2}} + {\alpha_{1}}\,\left( -3 + 4\,{\alpha_
{1}}\,\psi'({\alpha_{1}}) \right)  \right) \,
         \psi'({\alpha_{2}}) \right) }^2}\, .
 \end{eqnarray}

$\hfill \Box$
         \end{proposition}

\section{Applications}
The univariate gamma information geometry is known and
has been applied recently to represent and metrize
departures from randomness of, for example,
the processes that allocate gaps between occurrences of each amino acid
along a protein chain within the {\em Saccharomyces cerevisiae} genome,
see Cai et al~\cite{amino}, cosmological void distribution
and clustering of galaxies, and communications,
Dodson~\cite{ijtp,gsis,igcc}.

The new results on bivariate gamma geometry have
potential applications in any situation where positively correlated
variables $0<x<y<\infty$ are used to model a process with marginal
gamma distributions. In some of the applications mentioned above,
there are other variables associated and these may yield
some refinements of the existing models.
Two new case studies are being developed at present,
both involve bivariate data from measurements on stochastic porous media.

The first application concerns the structure of paper and nonwoven textiles,
in the manufacture of which fibres are deposited on a
continuous filter bed and form a near-planar bonded network. In the random case
it is easily seen that the mean number of sides is four for the polygonal voids
formed by the fibre process.
The distribution of polygonal void sizes is then given by the direct
product of independent identical exponential
distributions---reflecting the Poisson processes for fibre intersections.
For isotropic but {\em non-random} manufacturing processes,
Dodson and Sampson~\cite{aml} used the product of two
gamma distributions to obtain the void size distribution, recovering the
known random model as a special case. For some such isotropic materials it may be
appropriate to consider correlated polygon sides, since the voids tend to
be `roundish' suggesting positive correlation among polygon sides.
Some commercial processes involve preferential
alignment of fibres in the direction of manufacture, resulting in anisotropic
void distributions. In these cases the product of independent gamma
distributions needs to be replaced by the McKay bivariate model
and we shall report this study elsewhere.

The second application concerns tomographic images of soil samples in
hydrology surveys. Such images yield 3-dimensional reconstructions of
the porous structure and interestingly these generate a bivariate
positively correlated pair of random variables $0<x<y<\infty.$ Here,
the larger variable represents the size of a void and the smaller
variable represents the sizes of the throats that connect it to
neighbouring voids. This work is being pursued in collaboration with
Professor J. Scharkanski of the Federal University of Rio Grande do Sul, Brazil
and will be reported elsewhere~\cite{hydrology}.

The authors used {\em Mathematica} to perform analytic calculations and can make
available working notebooks for others to use.

\section{Concluding remarks}
We have formalised the concept of a process being `nearly random'
by proving that in the subspace information metric topology of Euclidean $\R^3,$
every neighbourhood of an exponential distribution contains a neighbourhood
of gamma distributions.
We have derived the information geometry of the 3-manifold $M$ of McKay
bivariate gamma distributions, which can provide a metrization of
departures from randomness and independence for bivariate processes.
Additionally, we give the metric for a 5-manifold version that
includes location parameters for the two random variables.
The curvature objects are derived for $M$ and those on three submanifolds
that illustrate some cases of possible practical interest.
As in the case of bivariate normal manifolds, we have negative scalar curvature
on the McKay bivariate gamma manifolds,
but here it is not constant and we show how it depends on correlation
in two cases.  These results have applications, for
example, in the characterization of stochastic materials.

\subsection*{Acknowledgement} The authors wish to thank Dr Hiroshi Matsuzoe
for helpful comments concerning affine immersions.

\end{document}